\documentclass[12pt]{article}
\usepackage{amsmath,amsfonts,latexsym,amssymb,amsthm,mathrsfs,array}

\usepackage{amsfonts}

\newcommand{\CC}{{\mathbb C}}
\newcommand{\ZZ}{{\mathbb Z}}
\newcommand{\NN}{{\mathbb N}}
\def\gg{{\mathfrak{g}}}
\def\hh{{\mathfrak{h}}}
\def\CC{{\mathbb C}}
\def\ZZ{{\mathbb Z}}
\newcommand{\fma}{\overset{\circ}{\mathfrak{g}}}
\newcommand{\m}{{\mathfrak{m}}}

\newcommand{\fmh}{\overset{\circ}{\mathfrak{h}}}
\newcommand{\np}{\overset{\circ}{n^{+}}}
\newcommand{\nm}{\overset{\circ}{n^{-}}}

\newcommand{\C}{{\mathbb C}}
\newtheorem{dfn}{Definition}[section]
\newcommand{\bdfn}{\begin{dfn}\rm}
\newcommand{\edfn}{\end{dfn}}
\newtheorem{thm}[dfn]{Theorem}
\newcommand{\bthm}{\begin{thm}}
\newcommand{\ethm}{\end{thm}}                   
\newtheorem{lmma}[dfn]{Lemma}                   
\newcommand{\blmma}{\begin{lmma}}                   
\newcommand{\elmma}{\end{lmma}}                   
\newtheorem{ppsn}[dfn]{Proposition}
\newcommand{\bppsn}{\begin{ppsn}}
\newcommand{\eppsn}{\end{ppsn}}
\newtheorem{crlre}[dfn]{Corollary}
\newtheorem{rmk}[dfn]{Remark}
\newcommand{\brmk}{\begin{rmk}} 
\newcommand{\ermk}{\end{rmk}}

\title{Weyl modules associated to Kac-Moody Lie algebras}
\author{ S. Eswara Rao, V. Futorny and Sachin S. Sharma}
\date{}
\begin{document}
\maketitle
\begin{abstract}
 Weyl modules were originally defined for affine Lie 
algebras by Chari and Pressley in \cite{CP}.
In this paper we extend the notion of Weyl modules for a Lie algebra
$\gg \otimes A$, where $\gg$ is any Kac-Moody algebra and A is any finitely
generated commutative associative algebra with unit over $\CC$,
and prove a tensor product decomposition theorem 
generalizing \cite{CP}.
\end{abstract}

\section{Introduction}
Let $\gg $ be a Kac-Moody Lie algebra and let $\hh$ be a Cartan subalgebra of $\gg$. 
Set $\gg' = [\gg,\gg]$ and $\hh' =\gg \cap \hh$. Let $\hh''$ be a vector subspace of 
$\hh$ such that $\hh' \oplus \hh''= \hh$. Let  
$A$ be a finitely generated commutative associative algebra with unit over $\CC$.
Denote $\stackrel{\sim}{\gg} \,= \gg' \otimes A \oplus \hh''$ and let
 $\gg = N^{-} \oplus \hh \oplus N^{+}$ be a standard triangular decomposition 
into positive and negative root subspaces and a Cartan subalgebra. 
 Let $\stackrel{\sim}{N}^- = N^{-} \otimes A, \,
\stackrel{\sim}{N}^+ = {N}^+  \otimes A$ and 
$\stackrel{\sim}{\hh} = \hh' \otimes A\oplus \hh''$.
Consider a linear map $\psi : \stackrel{\sim}{\hh} \rightarrow \C$.

   In \cite{CP} Chari and Pressley defined the Weyl modules for the loop 
algebras, which are nothing but the maximal integrable highest weight modules.
Feigin and Loktev \cite{FL} generalized the notion of Weyl module by replacing Laurent polynomial ring 
by any commutative associative algebra with unit and generalized the tensor decomposition theorem of \cite{CP}.
 Chari and Thang \cite{CTH} studied
Weyl modules for double affine Lie algebra. In \cite{CFK}, a functorial approach used
to study Weyl modules associated with the Lie algebra $\fma \otimes A$, where
$\fma$ is finite dimensional simple Lie algebra and A is a commutative algebra with unit over $\CC$.
Using this approach they \cite{CFK}
defined a Weyl functor from category  commutative associative algebra
 modules to simple Lie algebra modules, and studied tensor product properties of this functor. 
Neher and Savage \cite{ENSA} using generalized evaluation representation discussed more general case by 
replacing finite dimensional simple algebra with an infinite dimensional Lie algebra.

Let $\tau = \fma \otimes A_n \oplus \Omega_{A_n}/dA_n$ be a toroidal algebra, where
$\fma$ is finite dimensional simple Lie algebra and
$A_n = \CC[t_1^{\pm},\cdots, t_n^{\pm}]$ is a Laurent polynomial ring in commutating $n$ variables(see \cite{E}). 
It is proved
in \cite{ESTL} that any irreducible module with finite dimensional weight spaces of $\tau$ is in fact a module
for $\gg \otimes A_{n-1}$ where  $\gg$ is affinization of $\fma$. Thus it is important to study
$\gg \otimes A_{n-1}$-modules. Rao and Futorny \cite{EV} initiated the study of $\gg \otimes A_{n-1}$-modules
in their recent work. In our
paper we consider the $\gg \otimes A$-module, where $\gg$
is any Kac-Moody Lie algebra and $A$ is any finitely generated commutative associative algebra with unit over $\CC$.
Our work is kind of generalization of
the tensor product results in \cite{CFK, FL, CP}. 
For a cofinite ideal $I$ of $A$ we define a module  $M(\psi,I)$, and a Weyl module
$W(\psi, I)$ of $\stackrel{\sim}{\gg}$ (Section \ref{s2}). The main result of the paper is
the tensor product decomposition of $W(\psi, I)$, where $I$ is a finite intesection of
 maximal ideals.

The paper is organised as follows. We begin with preliminaries by stating some
basic facts about Kac-Moody algebras and Weyl modules.  In Section \ref{s1} 
 we define the modules $M(\psi,I)$ over $\stackrel{\sim}{\gg}$ and show that they have finite
 dimensional weight spaces and prove
  tensor decomposition theorem for them. Section \ref{s2} is devoted to the tensor decomposition
theorem for the Weyl module $W(\psi,I)$ over $\stackrel{\sim}{\gg}$.

\section{ PRELIMINARIES}
Let $\fma$ be a finite dimensional simple Lie algebra of rank $r$ with
a Cartan subalgebra $\fmh$. Let $ {\overset{\circ}{\Delta}}$ denote a root system of $\fma$ 
with respect to $\fmh$. Let ${\overset{\circ}{\Delta^+}}$ and ${\overset{\circ}{\Delta^-}}$ be a sets of
positive  and negative roots of $\fma$ respectively. 
Denote by $\alpha_1, \cdots , \alpha_r$ and $\alpha_1^{\vee},\cdots,\alpha_r^{\vee}$ 
 a sets of simple roots and simple coroots of $\fma$.
 Let $\fma = \np \oplus \fmh \oplus \nm$
be a triangular decomposition of $\fma$. Let $e_i$ and $f_i$ be the Chevalley generators of $\fma$. 
Let $\overset{\circ}{Q} = \oplus \, \ZZ \, \alpha_i$ and ${\overset{\circ}{P}} = \{\lambda \in \fmh^{\ast} : \lambda
(\alpha_{i}^{\vee})\in \ZZ \}$
be the root and weight lattice of $\fma$ respectively.
Set  ${\overset{\circ}{P_+}} = \{\lambda \in \fmh^{\ast} : \lambda(\alpha_{i}^{\vee})\geq 0 \}$, 
the set of dominant integral weights of $\fma$.
 
Recall that a $\fma$-module $V$ is said to be integrable if
it is $\fmh$-diagonalisable and all the Chevalley generators $e_i$ and 
$f_i$, $1 \leq i \leq r$, act locally nilpotently on $V$.
For commutative associative algebra with unit $A$, consider 
the Lie algebra algebra $\fma \otimes A$.
We recall the definition of  local Weyl module for $\fma \otimes A$ \cite{FL,CFK}.

\bdfn
Let $\psi : \fmh \otimes A \rightarrow \C$ be a linear map such that $\psi \mid_{ \fmh} = \lambda ,$ 
 $I$ a cofinite ideal  of $A$.
Then $W(\psi , I)$ is called  a local Weyl module for $\fma \otimes A$ if there
exists a nonzero $v \in W(\psi, I)$  such that\\
$U(\fma \otimes A)v = W(\psi, I), (\np \otimes A)v = 0, (h \otimes 1)v = \lambda(h)v$\\
$\psi\mid_{\fmh \otimes I} = 0,{ (f_i \otimes 1) ^{\lambda(\alpha_{i}^{\vee})+1}}v = 0$, for $i = 1, \cdots, r .$
\edfn

It is shown in \cite{FL}(Proposition 4) that the local Weyl modules exists and can be obtained
as quotient of global Weyl module.

Let  $\stackrel{\sim}{\gg}\, = \,\gg' \otimes A \oplus \hh'' $
is a Lie algebra with the following bracket operations:
\begin{align*}
[X \otimes a, Y\otimes b] &= [X,Y] \otimes ab ,\\
[h, X \otimes a] &= [h,X] \otimes a ,\\
[h,h'] &= 0 ,
\end{align*}
where $X, Y \in \gg'$, $h, h' \in \hh''$ and $a, b \in A$.
Let $\stackrel{\sim}{\hh} := \hh' \otimes A \oplus \hh''$ and $\stackrel{\sim}{\gg}\, = \,\stackrel{\sim}{N}^+
\oplus \stackrel{\sim}{\hh} \oplus  \stackrel{\sim}{N}^-$  be a triangular decomposition of $\stackrel{\sim}{\gg}$,
where $\stackrel{\sim}{N}^+ = {N}^+  \otimes A$ and $\stackrel{\sim}{N}^- = {N}^-  \otimes A$.

Let $\psi : \stackrel{\sim}{\hh} \rightarrow \CC$ be a linear map.
\bdfn
 A module $V$ of $\stackrel{\sim}{\gg}$ is 
called highest weight module (of highest weight $\psi$) if $V$ is generated by  a highest weight 
vector 
${v}$ 
such that 

\noindent 
(1) $\stackrel{\sim}{N}^+ v=0.$\\
(2) $h  \ {v}  =\psi (h) {v}$ for $h \in \ 
\stackrel{\sim}{\hh}, \psi \in \stackrel{\sim}{\hh}^{\ast}$.
\edfn

 Let $\CC$ be the one dimensional 
representation of $\stackrel{\sim}{N}^+ \oplus \stackrel{\sim}{\hh}$ where $\stackrel{\sim}{N}^+ $ acts 
trivially and $\stackrel{\sim}{\hh}$ acts via $h.1 = \psi(h) 1$
for $\forall h \in \stackrel{\sim}{\hh}$. Define the induced module 
$$\displaystyle{ M(\psi) = U(\stackrel{\sim}{\gg}) \displaystyle{\bigotimes_{U(\stackrel{\sim}{N}^+ \oplus 
\stackrel{\sim}{\hh})}} \CC} \,.$$
Then $M(\psi)$ is highest weight module and has a unique irreducible quotient
denoted by $V(\psi)$. 

\section{The modules $M(\psi, I)$ and its tensor decomposition} \label{s1}
Let $\alpha_1, 
\ldots, \alpha_l$ be a set of simple roots of $\gg$ and $\Delta^+$  a set of corresponding positive  roots.  
Let $Q = \displaystyle{\bigoplus_{i= 1}^{l}}{\ZZ \alpha_i}$ be root lattice of 
$\gg$ and  $Q_{+} = \displaystyle{\bigoplus_{i= 1}^{l}}{\ZZ_{\geq 0} \alpha_i}.$
Let $\lambda \in 
\hh^{\ast} $ be a dominant integral weight of $\gg$. Consider $\alpha \in \Delta^+$ and assume 
$\alpha =\sum n_i \alpha_i$.  
Define an usual ordering on $\Delta^+$ by $\alpha \leq \beta$ for $\alpha, 
\beta \in \Delta^+$ if $\beta - \alpha \in Q_{+}$.

 Let $I$ be a cofinite ideal of $A$. Let $\{ I_\alpha, \alpha \in \Delta^+\}$ be a sequence  of 
cofinite ideals of $A$ such that $I_{\alpha} \subseteq I$ and \\
(1) $\alpha \leq \beta \Rightarrow I_\beta \subseteq I_\alpha$.\\
(2) $I_\alpha I_\beta \subseteq I_{\alpha+\beta} $ if $\alpha +\beta \in 
\Delta^+$.

   For $\beta \in \Delta^+$ let $X_{-\beta}$ be a root vector 
corresponding to the root $-\beta$. For a cofinite ideal  $I$ of $A$ set 
$X_{-\beta} I = X_{-\beta} \otimes I$.

Let $\psi :\stackrel{\sim}{\hh} \rightarrow \C$ be a linear 
map such that $\psi\mid_{\hh' \otimes I} =0 $, $\psi\mid_{\hh}=\lambda \in 
\hh^*$ and $\lambda$ is  dominant integral.

\bdfn \label{df1}  We will denote by $M(\psi, \{I_\alpha, \alpha \in 
\Delta^+\})$  the  highest weight $\stackrel{\sim}{\gg}$-module with highest weight $\psi$ and highest weight
vector $v$  such that 
$(X_{-\beta} I_\beta)v=0$ for all $\beta \in \Delta^+$.
\edfn

We will show now the existence of  modules  $M(\psi, \{I_\alpha, \alpha \in 
\Delta^+\})$. Let $M(\psi)$ be 
the Verma module with a highest weight $\psi$ and a highest weight vector $v$. We will prove  that
the module generated by $X_{-\alpha}I_{\alpha}v$ is a proper submodule of $M(\psi)$ for all 
$\alpha \in \Delta_{+}$. We use induction on  the height of $\alpha$.
First recall that there is a cofinite 
ideal $I$ such that $\psi \mid_{\hh^{\prime} \otimes I} = 0$ and by definition 
$I_{\alpha} \subseteq I$ and $I_{\alpha} \subseteq I_{\beta}$ if $\beta \leq \alpha$.

Let us consider $X_{-\alpha_i} I_{\alpha_{i}}v$ for a simple root $\alpha_i$.
We will prove that $X_{-\alpha_i} a v$  generates a  proper submodule of $M(\psi)$
 where $a \in I_{\alpha_{i}}$.
 Indeed we have $X_{\alpha} b X_{-\alpha_i} a v = 0$ for any 
simple root $\alpha \neq \alpha_i$ and $b \in A$.  Let $N_i$ be the $\gg$-submodule generated by $X_{-\alpha_i}I_{\alpha_i}v$.
Then $X_{\alpha_i}b X_{-\alpha_i} a v = h_{i} ba v = 0$ as $I_{\alpha_{i}} \subseteq I$
and $\psi \mid_{\hh^{\prime} \otimes I} = 0$. So $M(\psi)/N_i \neq 0$ and the induction starts.

Let $\beta \in \Delta_+$ and $\mathrm{ht}(\beta) = n$.
Let $N $ be the  submodule generated by $\sum_{\nu \in \Delta_+}{X_{-\nu}I_{\nu}}v$ where $\mathrm{ht}{\nu} < n$.
Then by induction, $N$ is a proper submodule of $M(\psi)$.
Now consider $X_{\alpha_{i}}b X_{-\beta} a v$ where $\alpha_i$ is a simple root, $b \in A$ and $a \in I_{\beta}$.
But  $X_{\alpha_{i}}b X_{-\beta} a v = X_{-\beta + \alpha_i} ba v$. Since $\mathrm{ht}(\beta - \alpha_i) < n$
and $I_{\beta}$ is an ideal of $A$, we have $ba \in I_{\beta} \subseteq I_{\beta - \alpha_i}$. Hence,
we see that $X_{-\beta + \alpha_i} ba v \in N$. Therefore,
 $ X_{-\beta} a v$ is a
highest vector of $M(\psi)/N$, and hence generates a proper submodule.

\blmma \label{l1}
 Let $\gamma_1, \ldots, \gamma_n, \beta \in 
\Delta^+$.  Then $ B= X_{-\beta} I^{n+1}_\beta X_{-\gamma_1} a_1\ldots 
X_{-\gamma_n} a_n v=0$, for $ a_1, \ldots a_n  \in A$ and each $\gamma_i 
\leq \beta$.
\elmma

\begin{proof} 

We prove the statement by induction on $n$.  For $ n=0$ 
the lemma follows from the definition of the module.
We have $$B=X_{-\gamma_1} a_1 X_{-\beta} I^{n+1}_\beta X_{-\gamma_2} a_2 
\ldots X_{-\gamma_n} a_n {v} + [X_{-\beta}, X_{-\gamma_1}] 
I^{n+1}_\beta X_{-\gamma_2}a_2 \ldots X_{-\gamma_n} a_n {v}\,.$$
The first term is zero by induction on $n$.  Repeating the same argument 
$n$ times for the second term we end up with: \\

$B= [\ldots [X_{-\beta}, X_{-\gamma_{1}}], X_{-\gamma_2}], \ldots, X_{-\gamma_n}]I_\beta^{n+1} {v}$ .
Assume $$ [\ldots [X_{-\beta}, X_{-\gamma_{1}}], X_{-\gamma_2}], \ldots, X_{-\gamma_n}]\neq 0.$$ Then 
it is a nonzero multiple of
$
X_{-(\beta +\sum\gamma_i})$ and $
\beta +\sum\gamma_i$ is a root.
As each $\gamma_i \leq \beta$ we have $I_\beta \subseteq I_{\gamma_i}$. 

Thus
$$
I^{n+1}_\beta \subseteq I_\beta I_{\gamma_1} I_{\gamma_2}\ldots 
I_{\gamma_n} \subseteq I_{\beta +\sum \gamma_i}.
$$
Since
$$
X_{-(\beta +\sum\gamma_i}) I_{\beta +\sum \gamma_i} v=0,
$$
it completes the proof of the lemma.
\end{proof}

\bppsn 
 $M(\psi, \{I_\alpha, \alpha \in 
\Delta^+\})$ has finite dimensional weight spaces with respect to $\hh$.
\eppsn
\begin{proof}
  Follows from  Lemma \ref{l1} .
\end{proof}

 We now construct a special sequence of cofinite 
ideals $I_\alpha, \alpha \in \Delta^+$.  Let $I$ be any cofinite ideal 
of $A$.  Let $\psi: \stackrel{\sim}{\hh} \, \rightarrow \C$ be a 
linear map such that $\psi \mid_{\hh' \otimes I=0}$ and $\psi \mid_{\hh} 
=\gamma$ a dominant integral weight. Let recall that for $\alpha \in \Delta_{+}$ with $\alpha 
= \sum_{i = 1}^{l} {m_i \alpha_i}$, define $N_{\lambda, \alpha} = \sum_{i = 1}^{l}{m_i \lambda(\alpha_{i}^{\vee})}.$

Let $I_\alpha =I^{N_{\lambda, \alpha}}$. Now if
 $\alpha\leq\beta$ then it implies that , $N_{\lambda, \alpha} \leq N_{\lambda, \beta}$
 and hence $I_\beta \subseteq I_\alpha$. Suppose $\alpha, \beta \in 
\Delta^+$ such that $\alpha +\beta \in \Delta^+$.  Then clearly 
$I_\alpha I_\beta =I_{\alpha+\beta}$.  For this special sequence of ideals $I_{\alpha}$, define
$M(\psi, I):= M(\psi, \{I_\alpha, \alpha \in \Delta^+\})$.

 Let $I$ and $J$ be coprime cofinite ideals of 
$A$.
Consider linear maps $\psi_1, \psi_2:\stackrel{\sim}{\hh} \mapsto 
\C$  such that $\psi_1 \mid_{\hh'\otimes I} =0, \ \psi_2\mid_{\hh'\otimes J} = 0$, 
 $\psi \mid_{ \hh} =\lambda$ and $\psi_2 \mid_{\hh}=\mu$. Further assume that $\lambda$ 
and $\mu$ are dominant integral weights.  Let $M(\psi_1, I)$ and 
$M(\psi_2, J)$ be the corresponding highest weight modules. Now 
define the following  new sequence of cofinite ideals of $A$.
Let $K_\alpha =I^{N_{\lambda, \alpha}} \cap J^{N_{\mu, \alpha}} \subseteq 
I 
\cap J$.\\
It is easy to check that:\\
(1) If $\alpha\leq \beta, \ \alpha, \beta \in \Delta^+$ then $K_\beta 
\subseteq K_\alpha$.\\
(2) $K_\alpha \ K_\beta \subseteq K_{\alpha+\beta}$ if $\alpha, \beta, 
\alpha +\beta \in \Delta^+$.

Let $\psi =\psi_1 +\psi_2$, so that
$\psi\mid_{\hh' \otimes (I \cap J)} =0, \ \psi\mid_{\hh} =\lambda +\mu$. Then we have

\bthm \label{T1}
As a $\widetilde{\gg}$-module
$$ M(\psi_1 +\psi_2, \{K_\alpha, \alpha \in \Delta^+\}) \cong M(\psi, I) \otimes M(\psi_2, J) . $$
\ethm

The following is standard but we include the proof for convenience of the reader. 
\blmma \label{l2}
Let $I$ and $J$ be the coprime cofinite ideals of $A$. Then\\
a) $A=I^n +J^m$, for all $n,m \in \ZZ_{\geq 1}$.\\
b) $A/(I^n  \cap J^m)  
\cong$
$A/I^{n}  \oplus A/J^{m} .$

\elmma
\begin{proof}
 $a)$ As $I$ and $J$ are coprime there exist
$f \in I$ and $g \in J$ such that $f + g = 1$. Considering
the expression $(f+g)^{m+n+1} = 1$ , we see 
that the  left hand side is the sum of  two elements of 
$I^{m}$ and $I^{n}$. $b)$ is clear from $a)$ and the Chinese reminder theorem.
\end{proof}

Assume
$$
\begin{array}{lll}
\dim A/I^{N_{\lambda,\alpha}}  &=& m_\alpha\\
\dim A/J^{N_{\mu,\alpha}} &=&n_\alpha\\
\dim A/I^{N_{\lambda,\alpha}}  \cap J^{N_{\mu,\alpha}} &=&k_\alpha
\end{array}
$$
then 
$$
m_\alpha +n_\alpha =k_\alpha
$$
 by the above lemma.

\noindent {\bf Proof of  Theorem \ref{T1}.}

Let $a_{1,\alpha}, \ldots, a_{m_\alpha,\alpha}$ be a $\C$-basis 
of $A/I^{N_{\lambda,\alpha}}$.  Let $a^1_\alpha, a^2_{\alpha}, \ldots,$ 
be a $\C$-basis of $I^{N_{\lambda,\alpha}}$.  Then clearly 
$\stackrel{\sim}{N}^{-}$ has the following $\C$-basis:
$$
\begin{array}{l}
\{X_{-\alpha} a_{i,\alpha}, 1 \leq i \leq m_\alpha, \alpha \in 
\Delta^+\} \cup\\
\{ X_{-\alpha}  a^i_\alpha, i \in \NN, \ \alpha \in \Delta^+\}.
\end{array}
$$
Let $U_\lambda, \ U^\lambda$ be the subspaces of $U(\stackrel{\sim}{N}^-)$ spanned by the 
ordered products of the first set and the second set respectively.  Then by the 
PBW theorem we have 
$
U(\stackrel{\sim}{N}^-) = U_\lambda U^\lambda.
$

Let 
$$
\begin{array}{lll}
M &= &M(\psi_1 +\psi_2, \ \{K_\alpha, \alpha \in \Delta^+\}),\\
M_1 &=& M(\psi_1, I),\\
M_2&=& M(\psi_2, J).\\
\end{array}
$$
It is easy to see that 
$M_1=U_\lambda U^\lambda v =U_\lambda v $ as $U^\lambda v = \CC \,v$.  Since $M_1$ has finite dimensional weight 
subspaces we can define  the character of $M_1$ as follows:
$$
\mathrm{Ch}\,M_1 =\sum_{\eta\in Q^+} \dim M_{1, \ \lambda-\eta}\, e^{-(\lambda-\eta)}.
$$

Let $l_\alpha$ denote the multiplicity of the root $\alpha$.\\

It is standard that
\\
$\dim M_{1,\lambda-n} =K^1_\eta$, where $K^1_\eta$ is given by
$$
\prod_{\alpha \in \Delta^+}(1-e^{-\alpha)^{-m_\alpha l_\alpha}} 
=\sum_{\eta 
\in Q^+} K^1_\eta e^{-\eta}.
$$
Also we have that
$\dim M_{2,\mu-\eta} =K^2_\eta$, where
$$
\prod_{\alpha\in \Delta^+} (1-e^{-\alpha})^{-n_\alpha l_\alpha} 
=\sum_{\eta\in Q^+} K^2_\eta e^{-\eta},
$$
and 
$\dim M_{\lambda +\mu-\eta} =K_\eta$, where $K_\eta$ is given by 
$$
\prod_{\alpha \in \Delta^+}(1-e^{-\alpha})^{-k_\alpha \ l_\alpha} 
=\sum_{\eta\in Q^+} K_\eta e^{-\eta}
$$
(recall that $k_\alpha =m_\alpha +n_\alpha$).

From the above calculations we see that
$$
\dim \ M_{\lambda+\mu-\eta} =\dim(M_1\otimes M_2)_{\lambda +\mu-\eta}.
$$
Thus to prove the theorem it is sufficient to show that there is a 
surjective $\widetilde{\gg}$-homomorphism from $M$ to $M_1 \otimes M_2$.

Let $v_1$ and $v_2$ be the highest weight vectors of $M_1$ and $M_2$ 
respectively.  Let $U$ be  $\widetilde{\gg}$  submodule of $M_1 \otimes 
M_2$ generated by $v_1\otimes v_2$.  It is easy to check that $(\psi_1 
+\psi_2) (h' \otimes (I \cap J))=0$.  Recall that $K_\alpha 
=I^{N_{\lambda,\alpha}} \cap J^{N_{\mu,\alpha}}$.  We have 
$X_{-\alpha} K_\alpha (v_1\otimes v_2)=0$  which immediately implies that  $U$ is a quotient of $M$.
Hence to complete the proof of the theorem it is sufficient to prove that 
$ U \simeq M_1 \otimes M_2$.

Clearly, $M_1 \otimes M_2$ is linear span of vectors of the form $w_1 
\otimes w_2$ where
$$
\begin{array}{l}
w_1 = X_{-\lambda_1} a_1\ldots X_{-\lambda_n} a_n v,\\
w_2 =X_{-\beta_1} b_1\ldots X_{-\beta_m}b_m v_2 \,.
\end{array}
$$
Let $\beta\in \Delta^+$.  By the definition of $I_{\alpha}$ and the argument given in the proof of 
Lemma \ref{l1} it is easy to see that there 
exists $N >>0$ such that 
$$
X_{-\beta} I^N w_1=0, \ X_{-\beta}  J^N w_2 =0 .\leqno{(a)}
$$
Now recall that $A=I^N +J^N$. Let $1= f+g, \ f\in I^N, \ g\in J^N$.  For any
$h \in A$ write $ h= fh +gh$.  We will use induction on $m+n$.  Consider
$$
X_{-\beta} fh (w_1 \otimes w_2) =w_1 \otimes X_{-\beta} fh \ w_2 \ (by \ 
(a))
$$
$$
\begin{array}{l}
=w_1 \otimes X_{-\beta} (h-gh)w_2\\
=w_1 \otimes X_{-\beta} h \ w_2 \ by \ (by \ (a)).
\end{array}
$$
As $X_{-\beta} \ fh(w_1\otimes w_2) \in U$ (by induction $w_1 \otimes w_2 \in U$),
we conclude that $w_1\otimes X_{-\beta} h \ w_2 \in U$. Similarly we have $X_{-\beta} hw_1\otimes \ w_2 \in U$.
It easily follows now that
 $U\simeq M_1\otimes M_2$. This completes the proof of the 
theorem.

\section{Weyl modules for loop Kac-Moody algebras and its tensor decomposition}\label{s2}
  In this section we define maximal integrable highest 
weight modules for $\widetilde{\gg}$ and prove a  tensor product 
theorem for them.

Recall that $M(\psi, I)$ is a highest weight module with a highest weight 
$\psi$ and a highest weight vector $v$. Further $\psi\mid_{\hh' \otimes 
I}=0$ and $\psi\mid_{\hh} =\lambda$ a dominant integral weight. Let 
$\alpha_1,\ldots, \alpha_l$ be simple roots and $\alpha^{\vee}_1, \ldots 
\alpha^{\vee}_l$ be the simple coroots.  Let $\{ X_{\alpha_i}, \alpha^{\vee}_{i}, 
X_{-\alpha_i}\}$ be an $\mathfrak{sl}_2$ copy corresponding to the simple root 
$\alpha_i$.

\bdfn \label{dw}
Let $W$ be  
a highest weight $\widetilde{\gg}$-module with a highest weight $\psi$ 
and a highest weight vector $v$ such that \\
(1) $\psi\mid_{\hh' \otimes I}=0$,\\
(2) $\psi\mid_{\hh} =\lambda$,    \\
(3) $X_{-\alpha_i}^{\lambda(\alpha_i^{\vee})+1} v=0$ for $i=1,2, \ldots, l$.
\edfn

It follows immediately that $W$ is an integrable $\gg$-module (see \cite{K}). 
We will prove below that such module exists and 
has finite dimensional weight spaces.

By the result of \cite{FL} (see the proof of Proposition 6 and 16 of \cite{FL}) 
it follows that 
$$
X_{-\alpha_i} I^{\lambda(\alpha_i^{\vee})} \, v=0.
$$
Let $\alpha \in \Delta^+$ and $X_{-\alpha} =[X_{-\alpha_{i_1}}, [\ldots 
[X_{-\alpha_{i_{n-1}}}, X_{-\alpha_{i_n}}]]],$
where $\sum \alpha_{i_j} =\alpha$.
It is easy to check that $X_{-\alpha} I^{N_{\lambda,\alpha}} v=0$.  
Recall that $N_{\lambda, \alpha} =\sum n_i \lambda (\alpha_i^{\vee})$ if $\alpha =\sum n_i \alpha_i$.

Thus $W$ is an integrable quotient of $M(\psi, I)$. 
Denote by  $W(\psi, I)$ the 
maximal such quotient of $M(\psi, I)$ in the sense that  any 
integrable quotient of $M(\psi, I)$ is a quotient of $W(\psi, I)$.
In particular, $W(\psi, I)$ has finite dimensional weight spaces. We will 
prove at the end that for a cofinite ideal $I$ which is finite intersection of maximal ideals,
$W(\psi, I)$ is non-zero by explicitly  constructing 
its irreducible quotient. We will call $W(\psi, I)$ the \emph{Weyl module} associated with 
$\psi$ and $I$. Let $I$ and $J$ be coprime finite ideals of $A$.

\bthm \label{tw}
 $W(\psi_1 +\psi_2, \ I \cap J) \cong W(\psi, 
I) \otimes W(\psi, J)$ as $\widetilde{\gg}$-modules.
\ethm

To prove above theorem we need the following lemma.

\blmma $W(\psi_1, I) \otimes W(\psi_2, J)$ is a 
quotient of $W(\psi_1 +\psi_2, \ I \cap J)$.
\elmma
\begin{proof} 
 Let $v_1, v_2$ be  highest weight vectors of 
$W(\psi_1, I)$ and $W(\psi_2, J)$ respectively.
As in the earlier argument we can prove that $W(\psi_1, I) \otimes 
W(\psi_2, J)$ is a cyclic module generated by $v_1 \otimes v_2$. Recall 
that $W(\psi_1 +\psi_2, \  I \cap J)$ is a maximal integrable quotient 
of $M(\psi_1 + \psi_2, 
I \cap J)$.  But $W(\psi_1, I) \otimes W(\psi_2, J)$ is a integrable 
quotient of $M(\psi_1+ \psi_2, I \cap J)$.  Hence $W(\psi_1, I) \otimes 
W( \psi_2, J)$ is a quotient of $W(\psi_1 + \psi_2, I \cap J)$.
\end{proof} 

\emph{Proof of the Theorem \ref{tw}}: In view of the above lemma, 
 it is sufficient to prove that $W(\psi_1 + 
\psi_2, I \cap J)$ is a quotient of $W(\psi_1, I) \otimes W(\psi_2, J)$. 
Let $K_i$ be the kernel of the map $M(\psi_i, I) \to W(\psi_i, I)$, $i=1, 2$.  
Then it is a standard fact that $\tilde{K} = K_1 \otimes M(\psi_2, J) + 
M(\psi_1, I) \otimes K_2$ is the kernel of the map
$$M(\psi_1, I) \otimes M (\psi_2, J) \to W(\psi_1, I) \otimes W(\psi_2, J).$$
Let $V$ be any integrable quotient of $M(\psi_1, I) \otimes M(\psi_2, 
J)$ and $K$ be the kernel of the map $M(\psi_1, I) \otimes  M(\psi_2, 
J) \to V$.\\  
{\bf{Claim}} :   $\tilde{K} \ \subseteq K$.  This claim proves that $V$ is 
a quotient of $W( \psi_1, I) \otimes W(\psi_2, J)$.  In particular, 
$W(\psi_1 + \psi_2, I \cap J)$ is a quotient of $ W(\psi_1, I) \otimes 
W(\psi_2, J)$ which completes the proof of the theorem.\\ 
{\bf{Proof of the claim}} : Since $V$ is $\gg$-integrable, it follows 
that the set of weights of $V$ is $W$- invariant and it is contained 
in $\lambda + \mu - Q^+$.  Here $Q^+$ is a monoid generated by
 simple roots,  $W$ is  the Weyl  group corresponding 
to $\gg$.  Let $m = (\lambda + \mu )(\alpha_i^{\vee}) + 1$.  
Then $(X_{-\alpha_i} f)^m (v_1 \otimes v_2) = 0$ in $V$, $\forall \, f \in A$. Indeed, if this element is not zero then 
 $\lambda + \mu - m\alpha_i$ is a weight of $V$ implying that $\lambda + \mu + \alpha_i$ is also a weight
by the $W$-invariance property of weights. But this is a contradiction.
 Let now $N = \mathrm{max}\{(\lambda + \mu)(\alpha_i^{\vee})+1 : 1 \leq i \leq l\}$.
As $I^{N} + J^{N} = A$ by  Lemma \ref{l2}, choose $f \in I^{N}$ and $g \in J^{N} $ such that
$f + g = 1$. 

Consider
\begin{align*}
 B &=  (X_{- \alpha_i} f)^N (v_1 \otimes v_2)\ \\
& ={\displaystyle{\sum_{k_1 + k_2 = N}}} C_i(X_{-\alpha_i} 
f)^{k_1} v_1 \otimes (X_{-\alpha_i}f)^{k_2} v_2\in K,
\end{align*}
with some constants $C_i$.
Since $(X_{-\alpha_i} f) v_1 = 0$, it follows that $B = v_1 
\otimes (X_{-\alpha_i} f)^N v_2 \in K$ and
$B = v_1 \otimes (X_{-\alpha_i} (1 - g))^N v_2 \in K.$
But  $(X_{-\alpha_i} g) v_2 = 0$.  Hence 
$v_1 \otimes X^N_{-\alpha_i} v_2 \in K \,.$ 
Let $n_0$ be the least positive integer such that $v_1 \otimes 
X_{-\alpha_i}^{n_0} v_2 \in K$.  But then
$$X_{\alpha_i}(v_1 \otimes X_{-\alpha_i}^{n_0} v_2 ) =  v_1 \otimes X_{\alpha_i}X_{-\alpha_i}^{n_0} v_2
= (n_0(\gamma_i - n_0 + 1)) v_1 \otimes 
X_{-\alpha_i}^{n_0-1} v_2 \in K,$$
where $\gamma_i = \mu (\alpha_i^{\vee})$.  By the minimality of $n_0$ it 
follows that $\gamma_i + 1= n_0$. 
Thus we have proved that 
$v_1 \otimes X_{-\alpha_i}^{\mu(\alpha_i^{\vee})+1} v_2 \in K , \,\,\forall
\,\,i \,.$
Similarly we can  prove that
$X_{-\alpha_i}^{\lambda(\alpha_i^{\vee})+1} v_1 \otimes v_2 \in
K.$
Now by earlier argument we can conclude 
%(using the fact that $I$ and $J$ are relatively prime) 
that $\tilde{K} \subseteq K$. This completes the proof of the 
claim. 
 
\begin{crlre} \label{max} 
 Let $I$ be a cofinite ideal of $A$ such that $I = {\displaystyle{\bigcap_{i=1}^k}} \m_i$, where
$\m_i$'s for $1 \leq i \leq k$ are distinct maximal ideals of $A$.  
Let $\psi_1, \cdots \psi_k$ be linear maps from 
$\widetilde{\hh} \to \C$ such that $\psi_i\mid_{\hh' \otimes \m_i} = 0$
and $\psi_i\mid_{\hh} = \lambda_i$ 
a dominant integral weight.  Put $\psi = 
{\displaystyle{\sum_{k}}}\psi_i$ 
and 
$\lambda = \sum \lambda_i$. Then 
$$W(\psi, I) \cong {\displaystyle{\bigotimes_{i=1}^k}}W(\psi_i, 
\m_i).$$
\end{crlre}

 {\bf{Remark}}: See \cite{CFK}, for the similar tensor decomposition theorem of 
 the local Weyl module for $\fma \otimes A$, where A is any commutative associative algebra with unit. 
%{\bf{Fact}} : Suppose $\gg$ is simple finite dimensional Lie 
%Algebra.  Suppose $V$ is  a $\gg$-integrable module with finite 
%dimensional 
%weight spaces.  Then $V$ is necessarily finite dimensional.  

%To see this fact first observe that $V$ is completely reducible.  Also 
%note that the weights of each irreducible submodule is contained in a 
%coset of $P$ mod $Q$.  Each coset has a unique miniscule weight (see \cite{H})  
%which is a weight of the submodule.  There are only finitely many 
%cosets and hence finitely many miniscule weights.  If the $V$ decomposes 
%into 
%infinitely many component then the dimensions of the miniscule weights 
%will be infinite dimensional and that is a contradiction.

{\bf{Remark}}: Module $W(\psi, I)$ is  $\widehat{\gg}$-integrable.  
In fact, $(X^{\lambda(\alpha_i^{\vee})+1}_{-\alpha_i} \otimes f) v = 
0$ for $\alpha_i$ simple and $f \in A$.  Indeed, suppose 
 it is non-zero.  Then $\lambda - (\lambda (\alpha_i^{\vee}) + 
1)\alpha_i$ is a weight of $W(\psi, I)$.  Since the weight are 
$W$-invariant 
it follows  that $\lambda +\alpha_i$ is a weight which is impossible.

 We now construct irreducible quotients of  Weyl 
modules hence proving their existence. Let $I$ be a cofinite
ideal of $A$ such that $\displaystyle{I = \bigcap_{i = 1}^{p}{\m_{i}}}$, where 
$\m_i$ are distinct maximal ideals of $A$. Now as $A$ is finitely generated over $\CC$,
$A/\m_i \cong \CC$ for $1 \leq i \leq p$. So by Chinese reminder theorem, there is 
a surjective homomorphism from $A$ to $\displaystyle{\bigoplus_{i=1}^{p}{A/\m_i}}$. Hence we have
have a surjective homomorphism  $\Phi : \gg' \otimes A \rightarrow \displaystyle{\bigoplus_{i = 1}^{p}{\gg' \otimes A/\m_i}}
\cong \gg'_p = 
\gg' \oplus \cdots \oplus \gg'$(p-times) by $\Phi(x \otimes a) = (a_1 x,a_2 x, \cdots,a_p x)$, where
$(a_1,a_2, \cdots, a_p) \in \CC^{p}$ is a image of $a$ from the map $A \rightarrow \bigoplus_{i=1}^{p}{A/\m_i}
\rightarrow \CC^{p}$.  Now for $1 \leq i \leq p$, let $\psi_i$ be
a linear map from $\widetilde{\hh}$ to $\CC$ such that 
$\psi_i \mid \hh = \lambda_i$ where $\lambda_i \in P^{+}$. Then as 
${\displaystyle{\bigotimes_{i=1}^p}} V(\lambda_i)$ 
is an irreducible integrable module for $\gg_p'$, it is also irreducible $\gg' \otimes A$-module via $\Phi$
and so for $\stackrel{\sim}{\gg}$($\hh''$ acts on tensor product via $\psi$)
and the vector $v_{\lambda_1}\otimes \cdots \otimes v_{\lambda_p}$ and $\psi =
 \psi_1 + \cdots + \psi_p$ satisfy the the conditions
of the definition \ref{dw} with  $\displaystyle{I = \bigcap_{i = 1}^{p}{\m_{i}}}$.

{\bf {Open Problem}} :  Compute the character of $W(\psi, I)$ 
which is $W$-invariant. 

By  Corollary \ref{max} it is sufficient to compute the character of $W(\psi_i, 
\m_i$), where $\m_i$ is a maximal ideal.  

\bibliography{rnewweyl.bib}

\nocite*{}

\bibliographystyle{plain}

School of mathematics, Tata Institute of Fundamental Research,
Homi Bhabha Road, Mumbai 400005, India.\\
email: senapati@math.tifr.res.in

Instituto de Mathem\'atica e Estat\'\i stica, Universidade de S\~ao Paulo, S\~ao Paulo, Brasil.\\
email: futorny@ime.usp.br

School of mathematics, Tata Institute of Fundamental Research,
Homi Bhabha Road, Mumbai 400005, India.
email: sachin@math.tifr.res.in

\end{document}